\documentclass[12pt]{amsart}
\usepackage{amsmath, amssymb, amsthm}
\usepackage[dvips]{graphicx}
\theoremstyle{plain}
\newtheorem{thm}{Theorem}[]

\theoremstyle{definition}

\theoremstyle{remark}

\numberwithin{equation}{section}
\begin{document}
\baselineskip=19pt
\title{Realization of the mapping class group of handlebody 
by diffeomorphisms}
\author{Susumu Hirose}
\address{Department of Mathematics, 
Faculty of Science and Engineering, 
Saga University, 
Saga, 840--8502 Japan}
\email{hirose@ms.saga-u.ac.jp}
\thanks{This research was supported by Grant-in-Aid for 
Scientific Research (C) (No. 20540083), 
Japan Society for the Promotion of Science. } 
\begin{abstract}
For the oriented 3-dimensional handlebody constructed 
from a 3-ball by attaching $g$ 1-handles, 
it is shown that the natural surjection from the group of orientation preserving 
diffeomorphisms of it to the mapping class group of it 
has no section when $g$ is at least 6. 
\end{abstract}

\maketitle
%
Let $M$ be an $n$-dimensional compact oriented manifold and 
$S$ be a subset of $\partial M$. 
We denote the group of orientation preserving diffeomorphisms 
of $M$ whose restrictions on $S$ are identity by $\mathrm{Diff} (M, S)$, 
the subgroups of them consisting of elements that are isotopic to 
identity by $\mathrm{Diff}_0 (M, S)$, 
and the quotient group $\mathrm{Diff} (M, S)/\mathrm{Diff}_0 (M, S)$ 
by $\mathcal{M}(M,S)$. 
For an element $f$ of $\mathrm{Diff} (M, S)$, let $[f]$ be the element of 
$\mathcal{M}(M,S)$ represented by $f$. 
The homomorphism $\pi_{M,S}$ from $\mathrm{Diff} (M, S)$ to $\mathcal{M}(M,S)$ 
defined by $\pi_{M,S} (h) = [h]$ is a surjection. 
Let $\Gamma$ be a subgroup of $\mathcal{M}(M,S)$. 
We call a homomorphism $s$ from $\Gamma$ to $\mathrm{Diff} (M, S)$ 
which satisfies $\pi_{M,S} \circ s = id_{\Gamma}$ a {\em section\/} for 
$\pi_{M,S}$ over $\Gamma$. 
Morita \cite{Morita} showed that the natural surjection from 
$\mathrm{Diff}^2 (\Sigma_g)$ to the mapping class group $\mathcal{M}(\Sigma_g)$ of 
$\Sigma_g$ has no section over $\mathcal{M}(\Sigma_g)$ when $g \geq 5$. 
Markovic \cite{Markovic} (when $g \geq 6$) and 
Markovic and Saric \cite{MS} (when $g \geq 2$) 
showed that the natural surjection from 
$\mathrm{Homeo} (\Sigma_g)$ to $\mathcal{M}(\Sigma_g)$ has no section 
over $\mathcal{M}(\Sigma_g)$. 
By using the different method from them, 
Franks and Handel \cite{FK} showed that the natural surjection 
from $\mathrm{Diff} (\Sigma_g)$ to $\mathcal{M}(\Sigma_g)$ has no section 
over $\mathcal{M}(\Sigma_g)$ when $g \geq 3$. 

Let $H_g$ be an oriented 3-dimensional handlebody of genus $g$, 
which is an oriented 3-manifold constructed from a $3$-ball by attaching 
$g$ 1-handles. 
Let $\Sigma_g$ be an oriented closed surface of genus $g$, 
then $\partial H_g = \Sigma_g$. 
The restriction to the boundary defines a homomorphism 
$\rho_{\partial} : \mathrm{Diff}(H_g) \to \mathrm{Diff}(\Sigma_g)$, 
and $\rho_{\partial}$ induces a injection 
$\mathcal{M}(H_g) \hookrightarrow \mathcal{M}(\Sigma_g)$ 
since $H_g$ is an irreducible 3-manifold. 
We will show: 

\begin{thm}\label{thm:handle-lift}
If $g \geq 6$, 
there is no section for 
$\pi_{H_g} : \mathrm{Diff} (H_g) \to \mathcal{M}(H_g)$
over $\mathcal{M}(H_g)$. 
\end{thm}

For contradiction, we assume that there is a section 
$s : \mathcal{M}(H_g) \to \mathrm{Diff} (H_g)$. 
Let $\Gamma$ be a subgroup of $\mathcal{M}(H_g)$, and 
$i_{\Gamma}$ be the inclusion from $\Gamma$ to $\mathcal{M}(H_g)$. 
Then $\Gamma$ is a subgroup of $\mathcal{M}(\Sigma_g)$, 
and the composition $\rho_{\partial} \circ s \circ i_{\Gamma}$ is 
a section for $\pi_{\Sigma_g} : \mathrm{Diff}(\Sigma_g) \to \mathcal{M}(\Sigma_g)$ 
over $\Gamma$. 
Therefore, if we can find a subgroup $\Gamma$ of $\mathcal{M}(H_g)$, 
over which there is no section for $\pi_{\Sigma_g}$, 
then Theorem \ref{thm:handle-lift} follows. 

Let $D$ be a 2-disk in $\Sigma_g$, and $\Sigma_{g,1}$ be $\Sigma_g \setminus int \, D$. 
Let $c$ be an essential simple closed curve on $\Sigma_g$ such that 
$\Sigma_g \setminus c$ is not connected, 
then the closure of one component of $\Sigma_g \setminus c$ is 
diffeomorphic to $\Sigma_{g_1, 1}$ and 
the closure of the other component of $\Sigma_g \setminus c$ 
is diffeomorphic to $\Sigma_{g_2, 1}$. 
We remark that $g = g_1 + g_2$ and $g_1, g_2 \geq 1$. 
These diffeomorphisms induce injections 
$\mathcal{M}(\Sigma_{g_1, 1}, \partial \Sigma_{g_1, 1}) \to \mathcal{M}(\Sigma_g)$ 
and  $\mathcal{M}(\Sigma_{g_2, 1}, \partial \Sigma_{g_2, 1}) \to \mathcal{M}(\Sigma_g)$
(see \cite{PR}). 
By these injections, 
we consider $\mathcal{M}(\Sigma_{g_1, 1}, \partial \Sigma_{g_1, 1})$ and 
$\mathcal{M}(\Sigma_{g_2, 1}, \partial \Sigma_{g_1, 1})$ as subgroups of 
$\mathcal{M}(\Sigma_g)$. 
From Theorem 1.6 in \cite{FK} proved by Franks and Handel, we see: 

\begin{thm}\cite{FK}\label{thm:FK}
Let $\Gamma_1$ be a nontrivial finitely generated subgroup of 
$\mathcal{M}(\Sigma_{g_1,1}, \partial \Sigma_{g_1,1})$ 
such that $H^1 (\Gamma_1, \mathbb{R}) = 0$, 
and $\mu$ be an element of $\mathcal{M}(\Sigma_{g_2,1}, \partial \Sigma_{g_2,1})$ 
which is represented by a pseudo-Anosov homeomorphism on 
$int\, \Sigma_{g_2,1}$. 
Then there is no section for $\pi_{\Sigma_g} : \mathrm{Diff}(\Sigma_g) \to \mathcal{M}(\Sigma_g)$ over $\langle \Gamma_1, \mu \rangle$. 
\end{thm}

We assume $g \geq 6$. 
The 3-manifold $\Sigma_{2,1} \times [0,1]$ is diffeomorphic to $H_{4}$. 
Let $D_1$ be a 2-disk in $int\, \partial \Sigma_{2,1} \times [0,1] \subset 
\partial(\Sigma_{2,1} \times [0,1])$, $D_2$ and $D_3$ be 
disjoint 2-disks on $ \partial H_{g-6}$, and $D_4$ be 2-disk on $\partial H_2$. 
Along these 2-disks, we glue $\Sigma_{2,1} \times [0,1]$, $H_{g-6}$ and $H_2$ 
such that $D_1 = D_2$, $D_3=D_4$, 
then the $3$-manifold obtained as a result is diffeomorphic to $H_g$. 
By the above construction, we get two natural inclusions 
$\Sigma_{2,1} \times [0,1] \hookrightarrow H_g$ and 
$H_2 \hookrightarrow H_g$. 
These inclusions induce natural homomorphisms 
$i_1 : \mathcal{M}(\Sigma_{2,1} \times [0,1], \partial \Sigma_{2,1} \times [0,1]) 
\to \mathcal{M}(H_g) $ 
and $i_2 : \mathcal{M}(H_2, D_4) \to \mathcal{M}(H_g)$. 
If $[h]$ is in 
$\mathcal{M}(\Sigma_{2,1} \times [0,1], \partial \Sigma_{2,1} \times [0,1])$ 
(resp. $\mathcal{M}(H_2, D_4)$) 
represented by $h \in 
\mathrm{Diff}(\Sigma_{2,1} \times [0,1], \partial \Sigma_{2,1} \times [0,1])$ 
(resp. $\mathrm{Diff}(H_2,D_4)$), 
then $i_1([h])$ (resp. $i_2([h])$ is represented by extending $h$ to 
$H_g$ using the identity mapping 
on $H_g \setminus \Sigma_{2,1} \times [0,1]$ 
(resp. $H_g \setminus H_2$). 

We define homomorphisms 
$\Pi : \mathrm{Diff}(\Sigma_{2,1}, \partial \Sigma_{2,1}) \to 
\mathrm{Diff}(\Sigma_{2,1} \times [0,1], \partial \Sigma_{2,1} \times [0,1])$ 
by $\Pi(h) = h \times id_{[0,1]}$, 
and $I_1 : 
\mathrm{Diff}(\Sigma_{2,1} \times [0,1], \partial \Sigma_{2,1} \times [0,1]) 
\to \mathrm{Diff}(H_g)$ by the identity on 
$H_g \setminus \Sigma_{2,1} \times [0,1]$, 
then the composition $I_1 \circ \Pi$ induces an injection 
$P : \mathcal{M}(\Sigma_{2,1}, \partial \Sigma_{2,1}) \to 
\mathcal{M}(H_g)$. 
By applying Corollary 4.2 of \cite{PR} to the subsurface 
$\Sigma_{2,1} \times \{ 0,1 \} \subset \partial H_g$, the injectivity of $P$ is shown. 
Korkmaz \cite{Korkmaz} showed that 
$H_1(\mathcal{M}(\Sigma_{2,1}, \partial \Sigma_{2,1}), \mathbb{Z}) 
= \mathbb{Z}/10 \, \mathbb{Z}$, hence 
$H^1 (\mathcal{M}(\Sigma_{2,1}, \partial \Sigma_{2,1}), \mathbb{R}) 
=0$. Therefore, 
$\Gamma_1 = P(\mathcal{M}(\Sigma_{2,1}, \partial \Sigma_{2,1}))$ 
satisfies the assumption of Theorem \ref{thm:FK} when $g_1=g-2$, $g_2=2$. 

Fathi and Laudenbach \cite{FL} constructed a 
pseudo-Anosov homeomorphism $\phi$ on $\partial(H_2)$ which 
is a restriction of a homeomorphism on $H_2$. 
Definition of pseudo-Anosov homeomorphisms and 
terminologies (e.g., singular foliation) related to them 
can be found in \cite{CB}.  
Any pseudo-Anosov homeomorphism preserves the set of singular points 
of the singular foliation which is preserved by this homeomorphism. 
Since the number of singular points of singular foliation is finite, 
a proper power of $\phi$, say $\phi^n$, fixes some points. 
Let $p$ be a point fixed by $\phi^n$, then $\phi^n$ defines a 
pseudo-Anosov homeomorphism on 
$\partial(H_2) \setminus p = int\, \Sigma_{2,1}$. 
Let $\mu$ be an element of $\mathcal{M}(\Sigma_{2,1}, \partial \Sigma_{2,1})
\subset \mathcal{M}(\Sigma_g)$ 
represented by this homomorphism, then 
$\mu$ is an element of $\mathcal{M}(H_g)$ and satisfies the assumption of 
Theorem \ref{thm:FK} when $g_1=g-2$, $g_2=2$. 

Then $\langle P(\mathcal{M}(\Sigma_{2,1}, \partial \Sigma_{2,1})), \mu \rangle$ 
is a subgroup of $\mathcal{M}(H_g)$ and, 
by Theorem \ref{thm:FK}, there is no section 
$\langle P(\mathcal{M}(\Sigma_{2,1}, \partial \Sigma_{2,1})), \mu \rangle \to 
\mathcal{M}(\Sigma_g)$. 
Therefore, there is no section 
for $\pi_{H_g} : \mathrm{Diff} (H_g) \to \mathcal{M}(H_g)$ 
over $\mathcal{M}(H_g)$.


\end{document}